\documentclass[reqno,12pt,pdftex]{amsart}

\usepackage[top=25truemm,bottom=25truemm,left=25truemm,right=25truemm]{geometry}

\usepackage{amsmath,amsthm,amssymb,amsfonts,amscd,mathtools,color}
\usepackage[all]{xy} 
\usepackage[colorlinks,linkcolor=blue,citecolor=blue,urlcolor=red]{hyperref}
\numberwithin{equation}{section}

\theoremstyle{plain} 
	\newtheorem{thm}{Theorem}[section]
	\newtheorem*{thm*}{Theorem}
	\newtheorem{cor}[thm]{Corollary}
	\newtheorem{lem}[thm]{Lemma}
	
	\newtheorem{prop}[thm]{Proposition}

\theoremstyle{definition}
	\newtheorem{defn}[thm]{Definition}

\theoremstyle{remark}
	\newtheorem{rem}[thm]{Remark}
	
	\newtheorem*{pf}{Proof}
\def\CC{{\mathbb C}}

\def\HH{{\mathbb H}}

\def\PP{{\mathbb P}}
\def\QQ{{\mathbb Q}}
\def\RR{{\mathbb R}}

\def\ZZ{{\mathbb Z}}
\def\A{{\mathcal A}}

\def\D{{\mathcal D}}
\def\E{{\mathcal E}}

\def\O{{\mathcal O}}
\def\P{{\mathcal P}}

\def\S{{\mathcal S}}

\def\U{{\mathcal U}}

\def\Aut{{\rm Aut}}
\def\coh{{\rm coh}}

\def\Hom{{\rm Hom}}

\def\Stab{{\rm Stab}}
\def\gldim#1{{\rm gldim} \, #1}
\def\uSdim#1{\overline{\rm Sdim} \, #1}
\def\lSdim#1{\underline{\rm Sdim} \, #1}
\def\mod{{\rm mod}}

\def\ex#1{\langle #1 \rangle_{\rm ex}}

\begin{document}
\title{Global dimension of the derived category of an orbifold projective line}
\date{\today}
\author{Takumi Otani}
\address{Department of Mathematics, Graduate School of Science, Osaka University, Toyonaka Osaka, 560-0043, Japan}
\email{otani.takumi.sci@osaka-u.ac.jp}

\maketitle
\begin{abstract}
The global dimension of a triangulated category is defined to be the infimum value of the global dimensions of stability conditions on the triangulated category.
In this paper, we study the global dimension of the derived category of an orbifold projective line.
\end{abstract}
\section{Introduction}
A stability condition on a triangulated category is a generalization of the slope stability of coherent sheaves on a nonsingular projective curve and the King stability of modules over a finite-dimensional algebra.
Bridgeland introduced this notion and proved that the space of stability condition $\Stab(\D)$ is a complex manifold \cite{B}.
This complex manifold plays an important role in mirror symmetry.
In a suitable setting, it is expected that the complex manifold has a structure of a Frobenius manifold (cf.~\cite{Ta, BQS}).
It is a significant problem to construct invariants of Frobenius manifolds from triangulated categories.
Ikeda--Qiu introduced the {\em global dimension $\gldim{\sigma}$ of a stability condition $\sigma$} on a triangulated category, which is a generalization of the global dimension of a finite-dimensional algebra \cite{IQ}.
The {\em global dimension} $\gldim{\D}$ of a triangulated category $\D$ is defined to be the infimum value of the global dimension of a stability condition on $\D$, i.e., 
\begin{equation*}
\gldim{\D} \coloneqq \inf_{\sigma \in \Stab(\D)} \gldim{\sigma},
\end{equation*}
In a suitable setting, we expect that the global dimension of a triangulated category is equal to the conformal dimension of a Frobenius manifold (cf.~\cite{KOT, Q}).
In general, the global dimension is a non-negative real number.
Kikuta--Ouchi--Takahashi gave a classification of triangulated categories satisfying $\gldim{\D} < 1$.
The next problem is classifying triangulated categories whose global dimension is one.
There are several works related to the problem (cf.~\cite{KOT, Q, S}).
Kikuta--Ouchi--Takahashi also proved $\gldim{\D^b(C)} = 1$ for the derived category $\D^b(C)$ of a nonsingular projective curve $C$.
On the other hand, Qiu proved $\gldim{\D^b(Q)} = 1$ for the derived category $\D^b(Q)$ of a finite acyclic quiver $Q$.
In order to give a classification, we need to study more examples.
Motivated by Kikuta--Ouchi--Takahashi's work \cite{KOT}, we study the global dimension of the derived category of coherent sheaves on an orbifold projective lines as a next class of triangulated categories whose global dimension is one.
Let $A$ be a multiple of positive integers and $\Lambda$ a multiple of pairwise distinct elements of $\PP^1$.
One can define an orbifold projective line $\PP_{A, \Lambda}^1$ associated with $(A, \Lambda)$.
There are three types of orbifold projective lines, called domestic type, tubular type and wild type, which are determined by their Euler characteristic $\chi_A \in \QQ$.
In order to study the global dimension of $\D^b(\PP_{A, \Lambda}^1)$, we characterize (the $\widetilde{\rm GL}^+(2, \RR)$-orbit of) the slope stability from the viewpoint of stability conditions on the triangulated category.
\begin{thm}[{Theorem \ref{thm : main 1}}]\label{thm : main 1 in intro}
Let $\sigma$ be a stability condition on $\D^b(\PP_{A, \Lambda}^1)$. 
For $\tau \in \HH$, define a homomorphism $Z_\tau \colon K_0(\coh(\PP_{A, \Lambda}^1)) \longrightarrow \CC$ by 
\[
Z_\tau (E) = - \deg (E) + \tau \cdot {\rm rank} (E), \quad E \in K_0(\coh(\PP_{A, \Lambda}^1)).
\]
We have $\sigma = (Z_\tau, \coh(\PP_{A, \Lambda}^1))$ for some $\tau \in \HH$ if and only if $\sigma$ satisfies the following conditions:
\begin{enumerate}
\item $\O(\vec{x})$ for $0 \le \vec{x} \le \vec{c}$ are $\sigma$-semistable.
\item $S_\lambda$ for $\lambda \not\in \Lambda$ and $S_{i, j}$ for $i = 1, \dots, r$, $j \in \ZZ / a_i \ZZ$ are $\sigma$-semistable.
\item $\phi(S_\lambda) = 1,  \phi(S_{i,j}) = 1$ and $0 < \phi(\O) < 1$.
\item $m (S_\lambda) = a$ and $m (S_{i,j}) = a / a_i$, where $a = \mathrm{lcm} \, A$.
\end{enumerate}
\end{thm}
The global dimension of the above stability condition is important to compute $\gldim{\D^b(\PP_{A, \Lambda}^1)}$.
As an analogue of \cite[Theorem 5.16]{KOT}, we prove the following theorem.
\begin{thm}[{Theorem \ref{thm : main 2}}]
Let $\PP_{A, \Lambda}^1$ be an orbifold projective line of type $(A, \Lambda)$.
For $\tau \in \HH$, we denote by $\sigma_\tau$ the stability condition $(Z_\tau, \coh(\PP_{A, \Lambda}^1))$.
The following holds:
\begin{enumerate}
\item If $\chi_A > 0$, then the stability condition $\sigma_\tau$ on $\D^b(\PP_{A, \Lambda}^1)$ for $\tau \in \HH$ satisfies $\gldim{\sigma_\tau} = 1$.
\item If $\chi_A = 0$, then the stability condition $\sigma_\tau$ on $\D^b(\PP_{A, \Lambda}^1)$ for $\tau \in \HH$ satisfies $\gldim{\sigma_\tau} = 1$ and is of Gepner type with respect to $(\S, \uSdim{\D^b(\PP_{A, \Lambda}^1)})$.
\item If $\chi_A < 0$, then $\gldim{\sigma} > 1$ holds for any stability condition $\sigma$ on $\D^b(\PP_{A, \Lambda}^1)$.
\end{enumerate}
Moreover, for any orbifold projective line, it holds that 
\begin{equation*}
\gldim{\D^b(\PP_{A, \Lambda}^1)} = \uSdim{\D^b(\PP_{A, \Lambda}^1)} = 1,
\end{equation*}
where $\uSdim{\D^b(\PP_{A, \Lambda}^1)}$ is the upper Serre dimension of $\D^b(\PP_{A, \Lambda}^1)$.
\end{thm}
The conformal dimension of the Frobenius manifold constructed from the genus zero Gromov--Witten theory for $\PP_{A, \Lambda}^1$ is one, which coincides with $\gldim{\D^b(\PP_{A, \Lambda}^1)}$ as expected.
Recently, Qiu--Zhang \cite{QZ} studied stability conditions whose global dimension is one in the case $\chi_A > 0$.
Their result is a generalization of Theorem \ref{thm : main 1 in intro}.
\bigskip
\noindent
{\bf Acknowledgements.}
I am grateful to Atsushi Takahashi for his support and comments, and to Yu Qiu for informing us about his work \cite{QZ} and for pointing out a mistake in my previous paper.
I would like to thank Akishi Ikeda, Kohei Kikuta, Suiqi Lu, Yuuki Shiraishi and Hongxia Zhang for their fruitful discussions.
This work is supported by JSPS KAKENHI Grant Number JP21H04994.
\bigskip
\noindent
{\bf Notation.}
Throughout this paper, for a finite-dimensional $\CC$-algebra $A$ (resp. a quiver $Q$), the bounded derived category of finitely generated right $A$-modules (resp. $\CC Q$-modules) is denoted by $\D^b(A) \coloneqq \D^b \mod (A)$ (resp. $\D^b(Q) \coloneqq \D^b \mod (\CC Q)$).
We also denote the derived category $\D^b \coh(X)$ of coherent sheaves on an orbifold $X$ by $\D^b(X)$.
For a triangulated category $\D$, the group of autoequivalences of $\D$ is denoted by $\Aut(\D)$.
\section{Preliminaries}
\subsection{Stability condition}
Following \cite{B}, we recall stability conditions on a triangulated category.
Let $\D$ be a $\CC$-linear triangulated category of finite type.
Denote by $K_0(\D)$ the Grothendieck group of $\D$.
\begin{defn}[{\cite[Definition 1.1]{B}}]\label{defn : stability condition}
A {\em stability condition} on $\D$ consists of a group homomorphism $Z \colon K_0(\D)\longrightarrow\CC$, which is called the {\em central charge}, and the family of full additive subcategories $\P = \{ \P(\phi) \}_{\phi \in \RR}$, called the {\em slicing}, satisfying the following axioms:
\begin{enumerate}
\item if $E \in \P(\phi)$ then $Z(E) = m(E) \exp(\sqrt{-1} \pi \phi)$ for some $m(E) \in \RR_{>0}$.
\item for all $\phi \in \RR$, $\P (\phi + 1) = \P(\phi) [1]$,
\item if $\phi_1 > \phi_2$ and $A_i \in \P(\phi_i)$ then $\Hom_\D (A_1, A_2) = 0$, 
\item for each nonzero object $E \in \D$ there exists a finite sequence of real numbers
\[
\phi_1 > \phi_2 > \dots > \phi_n
\]
and a collection of triangles
\small
\[
\xymatrix{
0 = F_0 \ar[rr] & & F_1 \ar[rr] \ar[ld] & & F_2 \ar[r] \ar[ld] & \cdots \ar[r] & F_{n-1} \ar[rr] & & F_n = E \ar[ld] \\
& A_1 \ar@{-->}[ul] & & A_2 \ar@{-->}[ul] & & & & A_n \ar@{-->}[ul] & 
}
\]
\normalsize
with $A_i \in \P(\phi_i)$ for all $i = 1, \dots, n$.
\item (support property) There exists a constant $C > 0$ such that for all nonzero object $E \in \P(\phi)$ for some $\phi \in \RR$, we have 
\[
\| E \| < C | Z(E) |
\]
where $\| \cdot \|$ denotes a norm on $K_0(\D) \otimes_\ZZ \RR$.
\end{enumerate}
\end{defn}
Let $\sigma = (Z, \P)$ be a stability condition on $\D$.
A nonzero object $E \in \P(\phi)$ is called a {\em $\sigma$-semistable} of phase $\phi$, and a simple object of $E \in \P(\phi)$ is called a {\em $\sigma$-stable}.
For any interval $I \subset \RR$, we denote by $\P(I)$ the extension closure $\ex{\P(\phi) \mid \phi \in I}$.
The full subcategory $\P((0, 1])$ forms a heart in $\D$.
A stability condition can be described by the heart $\P((0, 1])$ and the central charge (see \cite[Proposition 5.3]{B}). 
We also denote a stability condition $\sigma = (Z, \P)$ with $\A = \P((0, 1])$ by $\sigma = (Z, \A)$.
Denote by $\Stab(\D)$ the space of all stability conditions on $\D$.
Bridgeland showed in \cite[Theorem 1.2]{B} that the space $\Stab(\D)$ is a complex manifold.
There are natural group actions on $\Stab(\D)$ commuting with each other \cite[Lemma 8.2]{B}.
The first one is a $\widetilde{\rm GL}^+(2, \RR)$-action, where $\widetilde{\rm GL}^+(2, \RR)$ is the universal covering space of ${\rm GL}^+(2, \RR)$. 
The $\widetilde{\rm GL}^+(2, \RR)$-action preserves semistable objects, but the phases are changed.
We can consider the $\CC$-action on $\Stab(\D)$ induced by the embedding $\CC \hookrightarrow \widetilde{\rm GL}^+(2, \RR)$.
The $\CC$-action is given by 
\begin{equation*}
s \cdot (Z, \P) = (e^{- \pi \sqrt{-1} s} \cdot Z, \P_{{\rm Re} (s)}), \quad s \in \CC, 
\end{equation*}
where $\P_{{\rm Re} (s)}(\phi) \coloneqq \P (\phi + {\rm Re} (s))$.
The other action is given by the group of autoequivalences $\Aut(\D)$ as follows:
\begin{equation*}
\Phi (Z, \P) = (Z \circ \Phi^{-1}, \Phi(\P)), \quad \Phi \in \Aut(\D).
\end{equation*}
\begin{defn}[{\cite[Definition 2.3]{T}}]
A stability condition $\sigma$ on $\D$ is called {\em of Gepner type} with respect to $(\Phi, s) \in \Aut(\D) \times \CC$ if the following condition holds:
\begin{equation*}
\Phi (\sigma) = s \cdot \sigma.
\end{equation*}
\end{defn}
\subsection{Global dimension of a triangulated category}
We recall the notion of a global dimension of a stability condition, which was defined by \cite{IQ}.
This is a generalization of the global dimension of a finite-dimensional algebra.
\begin{defn}[{\cite[Definition 5.4]{IQ}}]
For a stability condition $\sigma = (Z, \P)$ on $\D$, the {\em global dimension $\gldim{\sigma} \in \RR_{\ge 0} \cup \{ + \infty \}$ of $\sigma$} is given by
\begin{equation*}
\gldim{\sigma} \coloneqq \sup \{ \phi_2 - \phi_1 \in \RR \mid \Hom_\D(A_1, A_2) \ne 0 ~ \text{for} ~ A_i \in \P(\phi_i) \}.
\end{equation*}
\end{defn}
Note that the global dimension $\gldim{\sigma}$ of a stability condition $\sigma$ is preserved under the $\Aut(\D)$-action and $\CC$-action.
\begin{defn}[{\cite{Q}}]
Define the {\em global dimension $\gldim{\D} \in \RR_{\ge 0} \cup \{ + \infty \}$ of a triangulated category $\D$} by
\begin{equation*}
\gldim{\D} \coloneqq \inf_{\sigma \in \Stab(\D)} \gldim{\sigma}.
\end{equation*}
\end{defn}
The upper Serre dimension of a triangulated category $\D$ gives a lower bound of the global dimension of $\D$.
We refer the reader to \cite{KOT} for the Serre dimension.
\begin{prop}[{\cite[Theorem 4.2]{KOT}}]\label{prop : Serre dimension}
Let $\D$ be a triangulated category.
Assume that $\D$ is equivalent to the perfect derived category ${\rm perf}(A)$ of a smooth proper differential graded $\CC$-algebra.
Then, we have $\uSdim{\D} \le \gldim{\D}$.
\qed
\end{prop}
We recall an important result by Qiu.
\begin{prop}[{\cite[Proposition 3.5]{Q}}]\label{prop : totally semistable stability condition}
A stability condition $\sigma$ on $\D$ satisfies $\gldim{\sigma} \le 1$ if and only if every indecomposable object in $\D$ is $\sigma$-semistable.
\qed
\end{prop}
\section{Orbifold projective line}
In this section, we prove our main theorems.
We first recall the definition of an orbifold projective line.
We refer to \cite{GL} for the derived category of an orbifold projective line.
Let $r \ge 3$ be a positive integer.
Let $A = (a_1, \dots, a_r)$ be a multiplet of positive integers and 
$\Lambda = (\lambda_1, \dots, \lambda_r)$ a multiplet of pairwise distinct elements of $\PP^1$ normalized such that 
$\lambda_1 = \infty$, $\lambda_2 = 0$ and $\lambda_3 = 1$. 
Put $a \coloneqq {\rm lcm}(a_1, \dots, a_r)$.
\begin{defn}[{\cite[Section 1.1]{GL}}]
Let $r$, $A$ and $\Lambda$ be as above.
\begin{enumerate}
\item Denote by $L_A$ an abelian group generated by $r$-letters $\vec{x_i}$, $i = 1, \dots, r$ defined as the quotient 
\begin{equation*}
L_A \coloneqq \left. \bigoplus_{i = 1}^r \ZZ \vec{x}_i \middle/ \big\langle a_i \vec{x}_i - a_j \vec{x}_j \mid i, j = 1, \dots, r \big\rangle \right. .
\end{equation*}
\item Define an $L_A$-graded $\CC$-algebra $S_{A, \Lambda}$ by 
\begin{equation*} 
S_{A,\Lambda} \coloneqq \CC[X_1,\dots,X_r] / (X_i^{a_i} - X_2^{a_2} + \lambda_{i} X_1^{a_1} \mid i = 3, \dots, r.) ,
\end{equation*}
where $\deg (X_i) = \vec{x}_i$ for $i = 1, \dots, r$.
\end{enumerate}
\end{defn}
We call the element $\vec{c} \coloneqq a_1 \vec{x}_1 = \dots = a_r \vec{x}_r \in L_A$ the {\em canonical element} of $L_A$.
Since $L_A / \ZZ \vec{c} \cong \bigoplus_{i = 1}^r \ZZ/ a_i \ZZ$, each $\vec{x} \in L_A$ can be uniquely written in normal form
\[
\vec{x} = l \vec{c} + \sum_{i = 1}^r l_i \vec{x}_i, \quad l \in \ZZ, \quad 0 \le l_i < a_i.
\]
For $\vec{x} \in L_A$, we write $\vec{x} \ge 0$ when $l \ge 0$ and $l_i \ge 0$ for $i = 1, \dots, r$ in normal form.
We also write $\vec{x} \ge \vec{y}$ when $\vec{x} - \vec{y} \ge 0$.
\begin{defn}
Define a stack $\PP^1_{A,\Lambda}$ by
\begin{equation*}
\PP_{A, \Lambda}^1 \coloneqq \left[ \left( {\rm Spec} (S_{A, \Lambda}) \backslash \{ 0 \} \right) / {\rm Spec} ({\CC L_A}) \right],
\end{equation*}
which is called the {\em orbifold projective line} of type $(A, \Lambda)$. 
\end{defn}
Define a rational number $\chi_A \in \QQ$ by 
\begin{equation*}
\chi_A \coloneqq 2 + \sum_{i = 1}^r \left( \frac{1}{a_i} - 1 \right).
\end{equation*}
The rational number $\chi_A$ is called the {\em Euler characteristic} of $\PP_{A, \Lambda}^1$.
An orbifold projective line $\PP_{A, \Lambda}^1$ satisfying $\chi_A > 0$ (resp. $\chi_A = 0$ and $\chi_A < 0$) is called of {\em domestic type} (resp. {\em tubular type} and {\em wild type}).
Denote by $\mod^{L_A}(S_{A, \Lambda})$ the abelian category of finitely generated $L_A$-graded $S_{A, \Lambda}$-modules and denote by $\mod^{L_A}_0 (S_{A, \Lambda})$ the full subcategory of ${\rm gr}^{L_A}(S_{A, \Lambda})$ whose objects are $L_A$-graded finite length $S_{A, \Lambda}$-modules.
It is known by \cite[Section 1.8]{GL} that the abelian category $\coh(\PP^{1}_{A, \Lambda})$ of coherent sheaves on $\PP^{1}_{A, \Lambda}$ is given by 
\[
\coh(\PP^{1}_{A, \Lambda}) = \mod^{L_A} (S_{A, \Lambda}) / \mod^{L_A}_0 (S_{A, \Lambda}),
\]
The abelian category $\coh(\PP_{A, \Lambda}^1)$ is hereditary (see \cite[Section 2.2]{GL}).
Therefore, each indecomposable object in $\D^b(\PP_{A, \Lambda}^1)$ is given by a shift of an indecomposable object in $\coh(\PP_{A, \Lambda}^1)$.
Define a sheaf $\O (\vec{x})$ for $\vec{x} \in L_{A}$ by
\begin{equation*}
\O (\vec{x}) \coloneqq [ S_{A, \Lambda} (\vec{x}) ] \in \coh(\PP_{A, \Lambda}^1)
\end{equation*}
where $(S_{A, \Lambda}(\vec{x}))_{\vec{x}'} \coloneqq (S_{A, \Lambda})_{\vec{x} + \vec{x}'}$. 
Denote by ${\rm Pic} (\PP_{A, \Lambda}^1)$ the group consisting of (isomorphism classes of) sheaves $\O(\vec{x})$ for $\vec{x} \in L_A$ with multiplication induced by the tensor product.
It is known by \cite[Section 2.1]{GL} that there is an isomorphism of abelian groups
\[
L_A \cong {\rm Pic}(\PP_{A, \Lambda}^1), \quad \vec{x} \mapsto \O(\vec{x}).
\]
For each $\vec{x} \in L_A$, the sheaf $\O(\vec{x})$ induces an autoequivalence $- \otimes \O(\vec{x}) \in \Aut (\D^b(\PP_{A, \Lambda}^1))$.
\begin{prop}[{\cite[Section 2.2]{GL}}]\label{prop : Serre functor}
Define an element $\vec{\omega} \in L_A$ by 
\begin{equation*}
\vec{\omega} \coloneqq (r - 2) \vec{c} - \sum_{i =  1}^r \vec{x}_i.
\end{equation*}
Then, the autoequivalence $\S \coloneqq - \otimes \O (\vec{\omega}) [1] \in \Aut (\D^b (\PP_{A, \Lambda}^1))$ is the Serre functor.
\qed
\end{prop}
\begin{prop}[{\cite[Section 4.1]{GL}}]\label{prop : canonical tilting object}
The ordered set 
\[
\E \coloneqq ( \O, \O(\vec{x}_1), \cdots, \O((a_1 - 1)\vec{x}_1)), \dots, \O(\vec{x}_r), \dots, \O((a_r - 1)\vec{x}_r), \O(\vec{c}) )
\]
is a full strongly exceptional collection in $\D^b(\PP_{A, \Lambda}^1)$.
In particular, $T = \bigoplus_{0 \le \vec{x} \le \vec{c}} \O(\vec{x})$ is a tilting object.
\qed
\end{prop}
It follows from Proposition \ref{prop : canonical tilting object} that the set $\{ [\O(\vec{x})] \mid 0 \le \vec{x} \le \vec{c} \}$ is a basis of $K_0(\D^b(\PP_{A, \Lambda}^1))$.
\begin{defn}[{\cite[Section 2.5]{GL}}]
Let $\lambda \in \PP^{1} \setminus \Lambda$. 
Define $S_\lambda$ and $S_{i, j}$ for $i = 1, \dots, r$ and $j \in \ZZ / a_i \ZZ$ by the following exact sequences:
\begin{subequations}\label{eq : simple objects}
\begin{equation}\label{eq : skyscraper}
0\rightarrow \O \xrightarrow{X^{a_{1}}_{1} - \lambda X^{a_{2}}_{2}} \O (\vec{c}) \longrightarrow S_\lambda \rightarrow 0,
\end{equation}
\begin{equation}\label{eq : simple at an orbifold point}
0\rightarrow \O (j \vec{x}_{i}) \xrightarrow{X_{i}} \O ((j + 1) \vec{x}_{i}) \longrightarrow S_{i, j} \rightarrow 0.
\end{equation}
\end{subequations}
\end{defn}
%

%
The sheaf $\O(\vec{x})$ for $\vec{x} \in L_A$ is given explicitly in the Grothendieck group $K_0(\D^b(\PP_{A, \Lambda}^1))$.
\begin{lem}\label{lem : line bundle}
Let $\vec{x} = l \vec{c} + \sum_{i = 1}^r l_i \vec{x}_i \in L_A$ in normal form.
It holds that
\begin{eqnarray*}
[\O(\vec{x})] & = & [\O] + l [S_\lambda] + \sum_{i = 1}^r \sum_{j = 0}^{l_i - 1} [S_{i, j}], \\
\lbrack S_\lambda \rbrack & = & \sum_{j \in \ZZ / a_i \ZZ} [S_{i, j}], \quad i = 1, \dots, r,
\end{eqnarray*}
in $K_0(\D^b(\PP_{A, \Lambda}^1))$.
\end{lem}
\begin{pf}
Fix an element $\vec{x} = l \vec{c} + \sum_{i = 1}^r l_i \vec{x}_i \in L_A$.
By \cite[(2.5.3) and (2.5.4)]{GL}, for $\lambda \notin \Lambda$, $i = 1, \dots, r$ and $j \in \ZZ / a_i \ZZ$ we have 
\begin{equation}\label{eq : line bundle tensor}
S_\lambda \otimes \O(\vec{x}) \cong S_\lambda \quad \text{and} \quad S_{i, j} \otimes \O(\vec{x}) \cong S_{i, j + l_i},
\end{equation}
where we regard $j + l_i$ as an element of $\ZZ / a_i \ZZ$.
Since the exact triangle \eqref{eq : skyscraper} induces
\begin{subequations}\label{eq : ses for line bundle}
\begin{equation}
0 \to \O \Big( (k - 1) \vec{c} + \sum_{i = 1}^r l_i \vec{x}_i \Big) \longrightarrow \O \Big( k \vec{c} + \sum_{i = 1}^r l_i \vec{x}_i \Big) \longrightarrow S_\lambda \to 0
\end{equation}
for each $k \in \ZZ$, we have
\[
[\O(\vec{x})] = \Big[ \O \Big( \sum_{i = 1}^r l_i \vec{x}_i \Big) \Big] + l [S_\lambda].
\]
Since the exact triangle \eqref{eq : simple at an orbifold point} induces 
\begin{equation}
0 \to \O \Big( \sum_{k = 1}^{i - 1} l_k \vec{x}_k + j \vec{x}_i \Big) \longrightarrow \O \Big( \sum_{k = 1}^{i - 1} l_k \vec{x}_k + (j + 1) \vec{x}_i \Big) \longrightarrow S_{i, j} \to 0.
\end{equation}
\end{subequations}
for each $i = 1, \dots, r$ and $j \in \ZZ / a_i \ZZ$, we have 
\[
\Big[ \O \Big( \sum_{i = 1}^r l_i \vec{x}_i \Big) \Big] = [\O] + \sum_{i = 1}^r \sum_{j = 0}^{l_i} [S_{i, j}]
\]
Therefore, we obtain the statement.

The later statement follows from the short exact sequences \eqref{eq : skyscraper} and \eqref{eq : simple at an orbifold point}.
\qed
\end{pf}
The {\em rank} and {\em degree} are homomorphisms 
\[
{\rm rank} \colon K_0(\D^b(\PP_{A, \Lambda}^1)) \longrightarrow \ZZ, \quad {\rm deg} \colon K_0(\D^b(\PP_{A, \Lambda}^1)) \longrightarrow \ZZ
\]
defined as follows:
\[
{\rm rank}(\O) = 1, \quad {\rm rank}(S_\lambda) = 0, \quad {\rm rank}(S_{i,j}) = 0,
\]
\[
\deg(\O) = 0, \quad \deg(\O(j \vec{x}_i)) = \dfrac{a}{a_i} j, \quad \deg(\O(\vec{c})) = a.
\]
Note that the degree of the sheaf $\O(\vec{\omega})$ is given by $\deg(\O(\vec{\omega})) = - a \chi_A \in \ZZ$.
%

%
Define ${\rm vect}(\PP^{1}_{A, \Lambda})$ and $\coh_0 (\PP^{1}_{A, \Lambda})$ as the full subcategories of $\coh (\PP^{1}_{A, \Lambda})$ consisting of objects whose rank is greater than $0$ and finite length objects, respectively.
Each coherent sheaf $E \in \coh(\PP^{1}_{A, \Lambda})$ splits into a direct sum $E_+ \oplus E_0$ of $E_+ \in {\rm vect} (\PP^{1}_{A, \Lambda})$ and $E_0 \in \coh_0(\PP^{1}_{A, \Lambda})$.
The category $\coh_0 (\PP^{1}_{A, \Lambda})$ decomposes into a coproduct $\coprod_{\lambda \in \PP^1} \U_\lambda$, where $\U_\lambda$ denotes the uniserial category of finite length sheaves concentrated at the point $\lambda$ (\cite[Proposition 2.4 and 2.5]{GL}).
%

%
\subsection{Slope stability}
Geigle--Lenzing studied the {\em slope} stability as an analogue of the classical slope stability for a nonsingular projective curve \cite[Section 5]{GL}. 
Denote by $\HH$ the upper half-plane.
By the $\widetilde{\rm GL}^+(2, \RR)$-action on $\Stab(\D^b(\PP_{A, \Lambda}^1))$, we have the following
\begin{prop}[{cf.~\cite[Section 5]{GL}, \cite[Lemma 8.2]{B}}]\label{prop : slope stability}
For each $\tau \in \HH$, define a group homomorphism $Z_\tau \colon K_0(\D^b(\PP_{A, \Lambda}^1)) \longrightarrow \CC$ by 
\begin{equation}\label{eq : slope}
Z_\tau (E) = - \deg (E) + \tau \cdot {\rm rank} (E).
\end{equation}
Then, the pair $\sigma_\tau = (Z_\tau, \coh(\PP_{A, \Lambda}^1))$ is a stability condition on $\D^b(\PP_{A, \Lambda}^1)$.
In particular, an object $E \in \coh(\PP_{A, \Lambda}^1)$ is $\sigma_\tau$-semistable (resp. stable) if and only if $E$ is slope semistable (resp. stable).
\qed
\end{prop}
\begin{cor}[{cf.~\cite[Proposition 5.5]{GL}}]\label{cor : semistable objects}
Let $\sigma_\tau = (Z_\tau, \coh(\PP_{A, \Lambda}^1))$ be the stability condition in Proposition \ref{prop : slope stability}.
The sheaves $\O(\vec{x})$ for $0 \le \vec{x} \le \vec{c}$, $S_\lambda$ for $\lambda \not\in \Lambda$ and $S_{i, j}$ for $i = 1, \dots, r$, $j \in \ZZ / a_i \ZZ$ are $\sigma$-stable.

In particular, if $\chi_A \ge 0$, then $\gldim{\sigma_\tau} \le 1$.
\end{cor}
\begin{pf}
It is well-known that every line bundle is $\sigma_\tau$-stable (cf.~\cite[Section 5]{GL}).
Since $S_\lambda$ for $\lambda \not\in \Lambda$ and $S_{i, j}$ for $i = 1, \dots, r$, $j \in \ZZ / a_i \ZZ$ are simple objects in $\coh(\PP_{A,\Lambda}^1)$, they are $\sigma_\tau$-stable.

We see the later claim.
Recall that an indecomposable object is in ${\rm vect}(\PP_{A, \Lambda}^1)$ or $\coh_0 (\PP_{A, \Lambda}^1)$.
It follows from \cite[Proposition 5.5]{GL} that, when $\chi_A \ge 0$, indecomposable objects in ${\rm vect}(\PP_{A, \Lambda}^1)$ are $\sigma_\tau$-semistable.
By \eqref{eq : slope}, one can see that $\coh_0(\PP_{A, \Lambda}^1) = \P(1)$ (cf.~\cite[Section 5]{GL}).
In particular, any indecomposable objects in $\coh_0 (\PP_{A, \Lambda}^1)$ are $\sigma_\tau$-semistable.
Hence, the statement follows from Proposition \ref{prop : totally semistable stability condition}.
\qed
\end{pf}
Conversely, one can characterize stability conditions in Proposition \ref{prop : slope stability} from the viewpoint of stability conditions on a triangulated category.
\begin{prop}\label{prop : characterization of slope stabilities}
Let $\sigma = (Z, \P)$ be a stability condition on $\D^b(\PP_{A, \Lambda}^1)$.
Assume the following conditions: 
\begin{enumerate}
\item $\O(\vec{x})$ for $0 \le \vec{x} \le \vec{c}$ are $\sigma$-semistable.
\item $S_\lambda$ for $\lambda \not\in \Lambda$ and $S_{i, j}$ for $i = 1, \dots, r$, $j \in \ZZ / a_i \ZZ$ are $\sigma$-semistable.
\item $\phi(S_\lambda) = 1, \phi(S_{i,j}) = 1$ and $0 < \phi(\O) < 1$.
\end{enumerate}
Then, we have $\P(0, 1] \cong \coh(\PP_{A, \Lambda}^1)$.
\end{prop}
\begin{pf}
Let $\A = \P((0, 1])$.
Note that the heart $\A$ is closed under extensions.
For $\vec{y} \in L_A$, if we have $\O(\vec{y}) \in \A$, then \eqref{eq : simple objects} implies $\O(\vec{y} + \vec{c}), \O(\vec{y} + \vec{x}_i) \in \A$. 
In addition, for $\vec{y} \in L_A$, if we have $\O(\vec{y}) \in \A$ and $0 < \phi(\O(\vec{y})) < 1$, then it follows from $S_\lambda, S_{i, j} \in \P(1)$ and \eqref{eq : simple objects} that $\O(\vec{y} - \vec{c}), \O(\vec{y} - \vec{x}_i) \in \A$. 
Hence, by assumption, we obtain $\O(\vec{y}) \in \A$ for all $\vec{y} \in L_A$ inductively.
By \cite[Proposition 2.6]{GL}, each object $E \in {\rm vect}(\PP_{A, \Lambda}^1)$ has a filtration 
\[
0 = E_0 \subset E_1 \subset \dots \subset E_n = E,
\]
whose factors $E_k / E_{k - 1}$ are line bundles $\O(\vec{l}_k)$ for suitably chosen $\vec{l}_k \in L_A$.
Since $\O(\vec{x}) \in \A$ for all $\vec{x} \in L_A$ and $\A$ is extension-closed, we have $E \in \A$, hence ${\rm vect}(\PP_{A, \Lambda}^1) \subset \A$.
Each finite length object $F \in \U_{\lambda_i}$ for some $\lambda_i \in \Lambda$ (resp. $F \in \U_{\lambda}$ for some $\lambda \not\in \Lambda$) has a filtration
\[
0 = F_0 \subset F_1 \subset \dots \subset F_n = F,
\]
whose factors $F_k / F_{k - 1}$ are simple object $S_{i, j}$ for some $j \in \ZZ / a_i \ZZ$ (resp. simple object $S_\lambda$).
It follows from $S_{i, j}, S_\lambda \in \P(1)$ that $F \in \P(1)$, hence $\coh_0 (\PP_{A, \Lambda}^1) \subset \P(1)$.
Therefore, we have $\coh(\PP_{A, \Lambda}^1) \subset \A$.
Since $\coh(\PP_{A, \Lambda}^1)$ and $\A$ are hearts in $\D^b(\PP_{A, \Lambda}^1)$, we obtain $\coh(\PP_{A, \Lambda}^1) = \A$. 
\qed
\end{pf}
\begin{thm}\label{thm : main 1}
Let $\sigma = (Z, \P)$ be a stability condition on $\D^b(\PP_{A, \Lambda}^1)$.
We have $\sigma = \sigma_\tau$ for some $\tau \in \HH$ if and only if $\sigma$ satisfies the following conditions:
\begin{enumerate}
\item $\O(\vec{x})$ for $0 \le \vec{x} \le \vec{c}$ are $\sigma$-semistable.
\item $S_\lambda$ for $\lambda \not\in \Lambda$ and $S_{i, j}$ for $i = 1, \dots, r$, $j \in \ZZ / a_i \ZZ$ are $\sigma$-semistable.
\item $\phi(S_\lambda) = 1,  \phi(S_{i,j}) = 1$ and $0 < \phi(\O) < 1$.
\item $m (S_\lambda) = a$ and $m (S_{i,j}) = a / a_i$.
\end{enumerate}
\end{thm}
\begin{pf}
Assume that $\sigma$ satisfies the above four conditions.
Let $\tau = Z(\O) \in \HH$.
By Lemma \ref{lem : line bundle}, the fourth condition (4) implies 
\begin{eqnarray*}
Z(\O(j \vec{x}_i)) = Z(\O) - \dfrac{a}{a_i} j = Z_\tau(\O(j \vec{x}_i)), \quad 
Z(\O(\vec{c})) = Z(\O) - a = Z_\tau(\O(\vec{c})),
\end{eqnarray*}
which yields $Z = Z_\tau$.
Hence, the statement follows from Proposition \ref{prop : characterization of slope stabilities}.
\qed
\end{pf}
Recently, Qiu--Zhang gives a generalization of Proposition \ref{prop : characterization of slope stabilities} in \cite{QZ}.
They classified stability conditions whose global dimension is one in the case $\chi_A > 0$.
\subsection{Global dimension of the derived category}
We consider the global dimension of the derived category of an orbifold projective line.
It follows from Proposition \ref{prop : Serre functor} and \cite[Lemma 2.10]{DHKK} that 
\[
\uSdim{\D^b(\PP_{A, \Lambda}^1)} = \lSdim{\D^b(\PP_{A, \Lambda}^1)} = 1.
\]
The following theorem is an analogue of Kikuta--Ouchi--Takahashi's result \cite[Theorem 5.16]{KOT} for the case of curves.
\begin{thm}\label{thm : main 2}
Let $\PP_{A, \Lambda}^1$ be an orbifold projective line of type $(A, \Lambda)$.
The following holds:
\begin{enumerate}
\item If $\chi_A > 0$, then the stability condition $\sigma_\tau$ on $\D^b(\PP_{A, \Lambda}^1)$ for $\tau \in \HH$ satisfies $\gldim{\sigma_\tau} = 1$.
\item If $\chi_A = 0$, then the stability condition $\sigma_\tau$ on $\D^b(\PP_{A, \Lambda}^1)$ for $\tau \in \HH$ satisfies $\gldim{\sigma_\tau} = 1$ and is of Gepner type with respect to $(\S, \uSdim{\D^b(\PP_{A, \Lambda}^1)})$.
\item If $\chi_A < 0$, then $\gldim{\sigma} > 1$ holds for any stability condition $\sigma$ on $\D^b(\PP_{A, \Lambda}^1)$.
\end{enumerate}
Moreover, for any orbifold projective line $\PP_{A, \Lambda}^1$, it holds that 
\begin{equation*}
\gldim{\D^b(\PP_{A, \Lambda}^1)} = \uSdim{\D^b(\PP_{A, \Lambda}^1)} = 1.
\end{equation*}
\end{thm}
We prove Theorem \ref{thm : main 2} for the cases $\chi_A > 0$, $\chi_A = 0$ and $\chi_A < 0$ in Subsections \ref{subsec : domestic}, \ref{subsec : tubular} and \ref{subsec : wild}, respectively.
%

\subsubsection{Proof of Theorem \ref{thm : main 2} for the case $\chi_A > 0$}\label{subsec : domestic}
An orbifold projective line $\PP_{A, \Lambda}^1$ satisfies $\chi_A > 0$ if and only if $A = (1, p, q), (2, 2, r), (2, 3, 3), (2, 3, 4)$ or $(2, 3, 5)$, where $p, q \ge 1$ and $r \ge 2$.
By Corollary \ref{cor : semistable objects}, each indecomposable object in $\D^b(\PP_{A, \Lambda}^1)$ is $\sigma_\tau$-semistable.
Therefore, we obtain the statement for $\chi_A > 0$ by Proposition \ref{prop : Serre dimension} and \ref{prop : totally semistable stability condition}.
\begin{rem}
Geigle--Lenzing gave a triangle equivalence $\D^b(\PP_{A, \Lambda}^1) \cong \D^b(Q_A)$, where $Q_A$ is the extended Dynkin quiver $A_{p, q}^{(1)}$, $D_{r + 2}^{(1)}$, $E_6^{(1)}$, $E_7^{(1)}$ and $E_8^{(1)}$ corresponding to $A = (1, p, q)$, $(2, 2, r)$, $(2, 3, 3)$, $(2, 3, 4)$ and $(2, 3, 5)$, respectively (cf.~\cite[Section 5.4.1]{GL}).
Qiu proved in \cite[Theorem 5.2]{Q} that for an acyclic quiver $Q$, which is not of Dynkin type, there exists a stability condition $\sigma$ on $\D^b(Q)$ such that
\begin{equation*}
\gldim{\D^b(Q)} = \gldim{\sigma} = 1.
\end{equation*}
Hence, we can also obtain $\gldim{\D^b(\PP_{A, \Lambda}^1)} = 1$ by this argument.
\end{rem}

\subsubsection{Proof of Theorem \ref{thm : main 2} for the case $\chi_A = 0$}\label{subsec : tubular}
An orbifold projective line $\PP_{A, \Lambda}^1$ satisfies $\chi_A = 0$ if and only if $A = (2, 2, 2, 2), (3, 3, 3), (2, 4, 4)$ or $(2, 3, 6)$.
In particular, we have $a = 2, 3, 4$ and $6$ for $A = (2, 2, 2, 2), (3, 3, 3), (2, 4, 4)$ and $(2, 3, 6)$, respectively.
By Proposition \ref{prop : Serre functor}, the derived category $\D^b(\PP_{A, \Lambda}^1)$ is fractional Calabi--Yau such that $\S^a = [a]$ if and only if $\chi_A = 0$.
By Corollary \ref{cor : semistable objects}, each indecomposable object in $\D^b(\PP_{A, \Lambda}^1)$ is $\sigma_\tau$-semistable.
Hence, by Proposition \ref{prop : Serre dimension} and \ref{prop : totally semistable stability condition}, we have $\gldim{\sigma_\tau} = \uSdim{\D^b(\PP_{A, \Lambda}^1)} = 1$.

It follows from \cite[Theorem 5.6]{GL} and Proposition \ref{prop : Serre functor} that, for each $\phi \in \RR$, we have $\S(\P(\phi)) = \P(\phi + 1)$.
Hence we have $\S(\sigma_\tau) = 1 \cdot \sigma_\tau$, which means that $\sigma_\tau$ is of Gepner type with respect to $(\S, 1)$.
Therefore, we obtain the statement for $\chi_A = 0$.

\begin{rem}
Kikuta--Ouchi--Takahashi proved that for a fractional Calabi--Yau triangulated category $\D$, a stability condition $\sigma$ satisfies $\gldim{\sigma} = \uSdim{\D}$ if and only if $\sigma$ is of Gepner type with respect to $(\S, \uSdim{\D})$ \cite[Theorem 4.6]{KOT}.
The statement of Theorem \ref{thm : main 2} for $\chi_A = 0$ follows from this result.
\end{rem}

\subsubsection{Proof of Theorem \ref{thm : main 2} for the case $\chi_A < 0$}\label{subsec : wild}
\begin{lem}
For any stability condition $\sigma$ on $\D^b(\PP_{A, \Lambda}^1)$ it holds that $\gldim{\sigma} > 1$.
\end{lem}
\begin{pf}
Suppose that there exists a stability condition $\sigma$ such that $\gldim{\sigma} = 1$.
By Proposition \ref{prop : totally semistable stability condition}, every indecomposable object is $\sigma$-semistable.
Since it holds that $\Hom_{\D^b(\PP_{A, \Lambda}^1)} (\O ((n - 1) \vec{\omega}), \O(n \vec{\omega})[1]) \ne 0$ for any $n \in \ZZ$, we have 
\begin{equation*}
1 + \phi(\O(n \vec{\omega})) - \phi(\O((n-1) \vec{\omega})) \ge 1.
\end{equation*}
Repeating this argument shows that $\phi(\O(a \vec{\omega})) \ge \phi(\O)$.
On the other hand, since $- a \chi_A > 0$, we have $\Hom_{\D^b(\PP_{A, \Lambda}^1)} (\O, \O(- a \chi_A \vec{c})) \ne 0$, which gives $\phi (\O) \le \phi(\O(- a \chi_A \vec{c}))$.
Then, the equality $a \vec{\omega} = - a \chi_A \vec{c}$ yields a contradiction. 
\qed
\end{pf}
\begin{lem}
For any $\varepsilon > 0 $, there exists a stability condition $\sigma_\varepsilon$ on $\D^b(\PP_{A, \Lambda}^1)$ such that $\gldim{\sigma_\varepsilon} < 1 + \varepsilon$.
In particular, it holds that $\gldim{\D^b(\PP_{A, \Lambda}^1)} = 1$.
\end{lem}
\begin{pf}
By Proposition \ref{prop : slope stability} and $\deg(\O(\vec{\omega})) = - a \chi_A < 0$, the statement follows from the same argument of \cite[Lemma 5.15]{KOT}.
\qed
\end{pf}
Hence, we have finished the proof of Theorem \ref{thm : main 2}.
\qed

\end{document}